\documentclass[12pt]{article}
\usepackage[utf8]{inputenc}
\usepackage{graphicx, float, amsmath,amssymb,amsthm}
\usepackage{genyoungtabtikz, ytableau, multicol}
\usepackage{mathtools}
\usepackage{authblk, hyperref, fullpage, cite}
\definecolor{babyblue}{rgb}{0.54, 0.81, 0.94}
\definecolor{bittersweet}{rgb}{1.0, 0.44, 0.37}
\begin{document}

\theoremstyle{plain}
\newtheorem{theorem}{Theorem}
\newtheorem{corollary}[theorem]{Corollary}
\newtheorem{lemma}[theorem]{Lemma}
\newtheorem{proposition}[theorem]{Proposition}
\theoremstyle{definition}
\newtheorem{definition}[theorem]{Definition}
\newtheorem{example}[theorem]{Example}
\newtheorem{conjecture}[theorem]{Conjecture}
\newtheorem{remark}[theorem]{Remark}
\newtheorem*{thm}{Theorem}

\newcommand{\seqnum}[1]{\href{https://oeis.org/#1}{\rm \underline{#1}}}

\begin{center}
\vskip 1cm{\LARGE\bf 
The Gini index in the representation theory of the general linear group}
\vskip 1cm
\large
Grant Kopitzke \\
Department of Mathematics\\
         University of Wisconsin Stevens Point\\
         Stevens Point, WI 54481\\
	 USA \\
\href{mailto:gkopitzk@uwsp.edu}{\tt gkopitzk@uwsp.edu} \\
\end{center}

\vskip .2 in

\begin{abstract}
The Gini index is a function that attempts to measure the amount of inequality in the distribution of a finite resource throughout a population. It is commonly used in economics as a measure of inequality of income or wealth. We define a discrete Gini index on the set of integer partitions with at most $n$ parts and show how this function emerges in the representation theory of the complex general linear group.
 
\end{abstract}

%%%%%%%%%%%%%%%%%%%%%%%%%%%%%%%%%%%%%%%%%%%%%%%%%%%%%%%%%%%%%%%%%%%%%%%%%%
%Section 1 - Introduction
%%%%%%%%%%%%%%%%%%%%%%%%%%%%%%%%%%%%%%%%%%%%%%%%%%%%%%%%%%%%%%%%%%%%%%%%%%

\section{Introduction}
The Gini index, originally defined in 1912 by the Italian statistician Corrado Gini \cite{Gini:1}, is a function that measures statistical dispersion. Traditionally used in economics to measure how equitably a resource is distributed throughout a population, the Gini index has found a wide range of modern applications in fields including biology \cite{Biology:1}, the medical sciences \cite{Medical:1}, and climate science \cite{Climate:1}.

We will review the discretization of the Gini index to the set of partitions of a fixed positive integer $n$ (defined in \cite{Kopitzke:1}), and will extend this discretization to the set of partitions of a positive multiple of $n$ with $n$ parts. In section 3, we review the representation theory of the complex general linear group, $GL_n(\mathbb{C})$, of invertible $n\times n$ complex matrices -- with an emphasis on representations of $GL_n(\mathbb{C})$-harmonic polynomials. We will then show how the discrete Gini index appears in this setting as the top homogeneous degree in which an irreducible rational representation of the general linear group occurs as a factor in the $GL_n(\mathbb{C})$-harmonic polynomials (Theorem \ref{theorem1}). 

\subsubsection*{Acknowledgements}
We would like to thank the referee for their suggestions that significantly improved the content and exposition of this article.

%%%%%%%%%%%%%%%%%%%%%%%%%%%%%%%%%%%%%%%%%%%%%%%%%%%%%%%%%%%%%%%%%%%%%%%%%
%Section 2 - The Discrete Gini Index
%%%%%%%%%%%%%%%%%%%%%%%%%%%%%%%%%%%%%%%%%%%%%%%%%%%%%%%%%%%%%%%%%%%%%%%%%

\section{The discrete Gini index}

%Subsection 2.1 %%%%%%%%%%%%%%%%%%%%%%%%%%%%%

\subsection{The Gini index of an integer partition}

 The classical Gini index is usually defined in terms of Lorenz curves. The Lorenz curve of a distribution is the graph of a function $L$ for which $L(x)$ is the percentage of total wealth possessed by the poorest $x$ percent of the population \cite{Lorenz:1}. From this definition it follows that, regardless of the distribution, $L$ will be an increasing, convex function on the interval from $0$ to $1$, $L(0)=0$, and $L(1)=1$. 
 
 The data sets we use to measure inequality are necessarily discrete. To make the corresponding Lorenz curve continuous, the curve is often constructed by approximating the discrete data set with a smooth curve \cite{Jantzen:1}, or by linear interpolation \cite{Farris:1}. Bypassing this process, a discrete version of the Gini index was defined in \cite{Kopitzke:1} as follows. Let $\lambda=\left(\lambda_1\geq\lambda_2\geq\cdots\geq\lambda_n\right)$ be a partition of an integer $n$ -- that is, a finite tuple of $n$ non-negative decreasing integers whose sum is equal to $n$. Each such partition corresponds to a wealth distribution in which there are $n$ dollars distributed among $n$ people such that the most wealthy person has $\lambda_1$ dollars, and the least wealthy has $\lambda_n$ dollars. To measure the amount of inequality in the wealth distribution $\lambda$, the Gini index of $\lambda$ is then defined as 
\begin{equation}\label{eqn:gini}
    g(\lambda):=\binom{n+1}{2}-\sum_{i=1}^ni\lambda_i.
\end{equation}

This definition is best understood in terms of discrete Lorenz curves. The discrete Lorenz curve of an integer partition $\lambda$ of $n$ is defined on the interval from 0 to $n$ by
\begin{equation}
    L_{\lambda}(x):= \begin{cases} 
      0 & x= 0 \\
      \displaystyle{\sum_{i=0}^{j-1}}\lambda_{n-i} & x\in (j-1,j] \\
   \end{cases}
\end{equation} 
for $j=1,\ldots,n$. 
\begin{example}
\label{example1}
Suppose there are $5$ dollars distributed among a population of $5$ people according to the partition $(3 \geq 1 \geq 1 \geq 0 \geq 0)$. The most equitable distribution would be the so-called \emph{flat} partition $(1 \geq 1 \geq 1 \geq 1 \geq 1)$ in which each person has the same amount of money; the Lorenz curve of this distribution is sometimes called the \emph{line of equity} or the \emph{line of equality}. The Lorenz curves of each of these distributions are pictured in Figure~\ref{figure1}.\\
\begin{figure}[H]
    \centering
    \includegraphics[width=15cm]{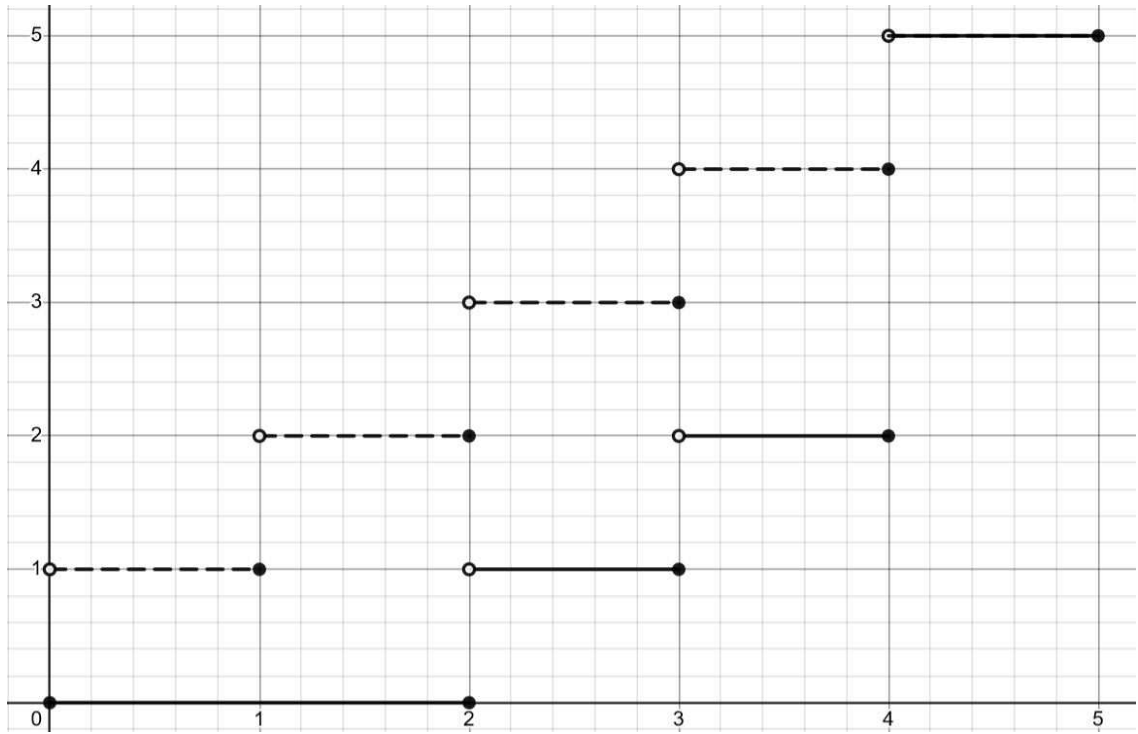}
    \caption{The line of equality (dashed) and the Lorenz curve of the partition ${(3 \geq 1 \geq 1 \geq 0 \geq 0)}$ of $5$ (solid).}
    \label{figure1}
\end{figure}
\end{example}
The Gini index of a partition $\lambda$ can easily be seen as the difference in the area between the line of equality and the curve $L_{\lambda}(x)$. Indeed, one computes that the area under the line of equality is $\binom{n+1}{2}$, and the area under $L_{\lambda}(x)$ is $\sum_{i=1}^n i\lambda_i$. In Example~\ref{example1} we see that the Gini index of $(3 \geq 1 \geq 1 \geq 0 \geq 0)$ is 7. The Lorenz curve of $\lambda$ will always coincide with the line of equality on the interval $(n-1,n]$, so the formula in Equation \ref{eqn:gini} simplifies to:
\begin{equation}\label{eqn:gini1}
    g(\lambda)=\binom{n}{2}-\sum_{i=1}^n(i-1)\lambda_i.
\end{equation}  
The function $\sum_{i=1}^n(i-1)\lambda_i$ appears frequently throughout algebraic combinatorics and is sometimes referred to as the weighted total of the partition $\lambda$ \cite{Macdonald:1}; thus we will adopt the function notation
\begin{equation}\label{eqn:weightedsum}
    b(\lambda):=\sum_{i=1}^n(i-1)\lambda_i.
\end{equation}

%Subsection 2.2 %%%%%%%%%%%%%%%%%%%%%%%%%%%%%%

\subsection{A generalization of the discrete Gini index}
The definition of the Gini index in Equation \ref{eqn:gini1} is unnecessarily restrictive, and can easily be generalized to allow for dollar amounts that are multiples of $n$. To that end, let $n$ and $k$ be positive integers, and consider a wealth distribution in which $nk$ dollars are distributed among $n$ people. The possible distributions are in one-to-one correspondence with the partitions of $nk$ with at most $n$ parts. 

The discrete Gini index can be generalized to this setting by extending the notion of the discrete Lorenz curve, and defining the Gini index of a distribution to be the area between the line of equality and the discrete Lorenz curve of the distribution. The Gini index in this setting is best understood via the following example. 

\begin{example}
\label{example2}
Suppose $15$ dollars are distributed among $5$ people according to the partition $\lambda=(6 \geq 4 \geq 3 \geq 1 \geq 1)$. The line of equality is the discrete Lorenz curve of the most equitable distribution, which in this setting is the flat partition $(3 \geq 3 \geq 3 \geq 3 \geq 3)$, or $(3^5)$ in simplified form. The discrete Lorenz curves for these distributions are pictured in Figure~\ref{figure2}.
\begin{figure}[H]
    \centering
    \includegraphics[width=15cm]{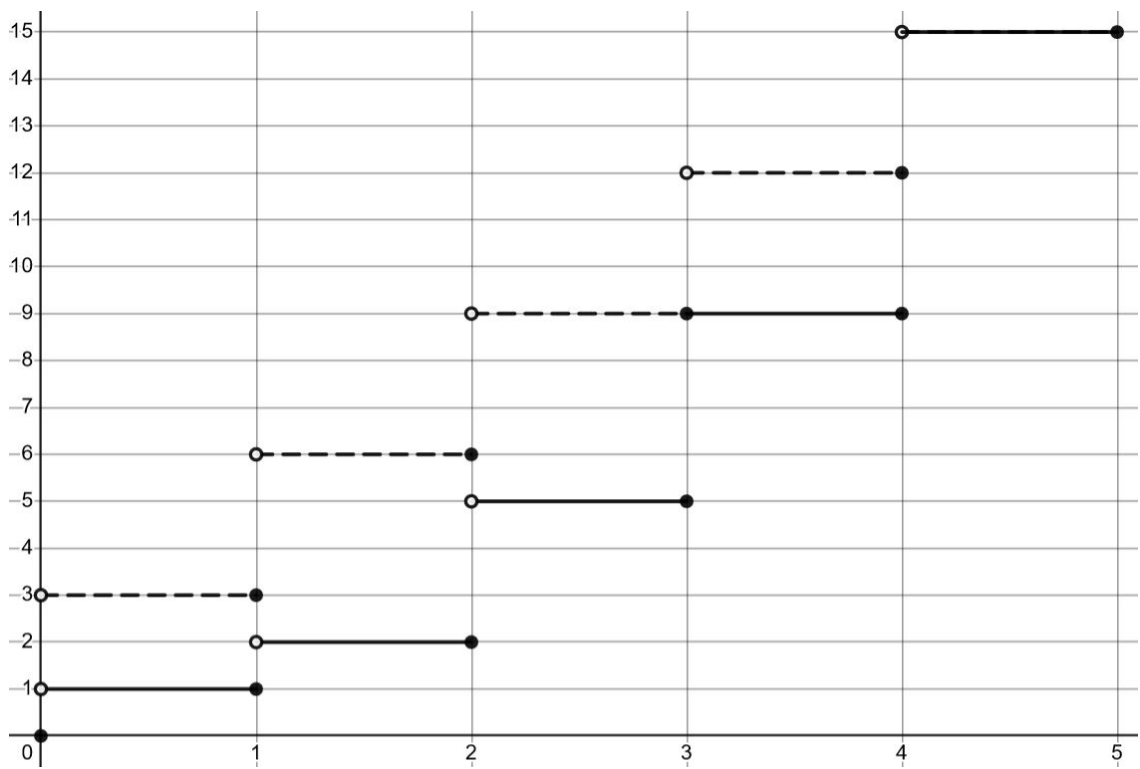}
    \caption{The line of equality (dashed) and the discrete Lorenz curve (solid) of the partition $(6 \geq 4 \geq 3 \geq 1 \geq 1)$ of $15$ with $5$ parts.}
    \label{figure2}
\end{figure}
\end{example}
We define the Gini index $g_{15,5}(\lambda)$ of this distribution to be the area between the line of equality and the discrete Lorenz curve of $\lambda=(6 \geq 4 \geq 3 \geq 1 \geq 1)$: 
\[b\left((3^5)\right)-b\left(\lambda\right)=13.  \]
In general, if $n$ and $k$ are positive integers and $\lambda$ is a partition of $nk$ with at most $n$ parts, then the Gini index of $\lambda$ is the area between the line of equality (the discrete Lorenz curve of the flat partition $(k^n)$) and the discrete Lorenz curve of $\lambda$:
\begin{equation}\label{eqn:gini2}
    g_{nk,n}(\lambda):=b\left((k^n)\right)-b(\lambda).
\end{equation} 
Notice that in the case $k=1$, the function in Equation \ref{eqn:gini2} specializes to the Gini index defined in Equation \ref{eqn:gini1}.

%%%%%%%%%%%%%%%%%%%%%%%%%%%%%%%%%%%%%%%%%%%%%%%%%%%%%%%%%%%%%%%%%%%%%%%%%
%Section 3 - Representation Theory
%%%%%%%%%%%%%%%%%%%%%%%%%%%%%%%%%%%%%%%%%%%%%%%%%%%%%%%%%%%%%%%%%%%%%%%%%

\section{Representation theory}
This section provides a brief overview of the basic representation theory of the complex general linear group and the structures relevant to our main result. Even so, some familiarity with the representation theory of classical groups is assumed. 

% Subsection 3.1 - Representation theory of GLN %%%%%%%%%%%%%%%%%%%%%%%%%

\subsection{Representation theory of $GL_n(\mathbb{C})$}

%Define the basics for rep. theory. of GL_n ala Fulton, then provide Fulton's theorem of highest weights.

%Precise statement of the theorem of the highest weight
The complex general linear group, $GL_n(\mathbb{C})$, is the set of all $n\times n$ complex invertible matrices, equipped with the operation of matrix multiplication. A finite dimensional \footnote{Unless otherwise stated, all representations in this paper are finite dimensional.} representation of $GL_n(\mathbb{C})$ is a homomorphism $\rho:GL_n(\mathbb{C})\longrightarrow GL(V)$, where $GL(V)$ is the group of invertible linear transformations of a finite dimensional complex vector space $V$. The representation is \emph{rational} if, after choosing a basis for $V$, $\rho(g)$ is a matrix whose entries are fixed rational functions in the matrix entries of $g$, where $g\in GL_n(\mathbb{C})$. When the action of $\rho$ is understood, the vector space $V$ itself is often referred to as the representation. 

An irreducible representation (irrep) is a nonzero representation with no proper nontrivial subrepresentations. In other words, $V$ is irreducible if there are no non-trivial subspaces $W\subset V$ that are closed under the action of $GL_n(\mathbb{C})$ restricted to $W$. The general linear group is \emph{reductive}, which implies that all rational representations can be written as a direct sum of irreducible rational representations. Irreps of $GL_n(\mathbb{C})$ are classified by the Theorem of the Highest Weight, which we now review.

Let $H$ denote the subgroup of $GL_n(\mathbb{C})$ consisting of diagonal matrices. Let $V$ be a rational representation of $GL_n(\mathbb{C})$. A vector $v\in V$ has \emph{weight} $\alpha=(\alpha_1,\ldots,\alpha_n)$ in $\mathbb{Z}^n$ if 
$\rho(x)(v)=x_1^{\alpha_1}\cdots x_n^{\alpha_n}v$
for all $x=\text{diag}(x_1,\ldots,x_n)$ in $H$.

Let $B$ denote the Borel subgroup of upper triangular matrices in $GL_n(\mathbb{C})$. The weight $\alpha$ of a vector $v\in V$ is called a \emph{highest weight} of $V$ if $\rho(B)(v)=\rho(\mathbb{C}^*)(v)$, where $\mathbb{C}^*$ is the set of all non-zero complex numbers. Two irreducible rational representations are isomorphic if and only if they have the same highest weight. These concepts are summarized in the Theorem of the Highest Weight.

\begin{thm}[Theorem of the Highest Weight]
For every tuple $\alpha=(\alpha_1,\ldots,\alpha_n)\in\mathbb{Z}^n$ with $\alpha_1\geq \ldots\geq \alpha_n$, there is a unique irreducible rational representation $V^{\alpha}$ of $GL_n(\mathbb{C})$ with highest weight $\alpha$. Moreover, all such irreducible rational representations of $GL_n(\mathbb{C})$ are of this form.
\end{thm}

For more on the representation theory of $GL_n(\mathbb{C})$ we refer the reader to \cite{Fulton:1}, \cite{Stanley:1}, and \cite{Wallach:1}.

% Subsection 3.2 - Harmonic Polynomials %%%%%%%%%%%%%%%%%%%%%%%%%%%%%%

\subsection{$GL_n(\mathbb{C})$--Harmonic polynomials}

%CHECK THIS CAREFULLY - a lot of it was summarized from Jeb's work, so check for plagiarism

%The harmonic analysis of GLn, ending with a formula for the graded multiplicity using Kostant's partition function

The discrete Gini index arises as the top degree in which an irrep appears as a direct summand in the space of harmonic polynomials on the Lie algebra of $GL_n(\mathbb{C})$. In this section, we provide a cursory overview of the $GL_n(\mathbb{C})-$harmonic polynomials.

It is well-known that the Lie algebra $\mathfrak{g}$ of $G=GL_n(\mathbb{C})$ is the Lie algebra of $n\times n$ complex matrices, $M_n(\mathbb{C})$. The general linear group acts on its Lie algebra by the adjoint representation,
\[ Ad(g)X:=gXg^{-1}, \]
for $g\in G$ and $X\in\mathfrak{g}$. Choose a basis $X_1,\ldots X_{n^2}$ for $\mathfrak{g}$, and define the algebra of polynomial functions on $\mathfrak{g}$ by identifying 
\[\mathbb{C}[\mathfrak{g}]=\mathbb{C}[X_1,\ldots,X_{n^2}].\]
The group $G$ acts on $\mathfrak{g}$ via the adjoint representation, and therefore also acts on $\mathbb{C}[\mathfrak{g}]$ by conjugation. The algebra $\mathbb{C}[\mathfrak{g}]$ is thus an infinite dimensional graded representation of $G$ with gradation
\[\mathbb{C}[\mathfrak{g}]=\bigoplus_{d\geq 0}\mathbb{C}[\mathfrak{g}]_d,\]
where $\mathbb{C}[\mathfrak{g}]_d$ is the vector space of homogeneous degree $d$ polynomials in $\mathbb{C}[\mathfrak{g}]$. The ring of $G-$invariant polynomials in $\mathbb{C}[\mathfrak{g}]$ is the set
\[ \mathbb{C}[\mathfrak{g}]^G=\left\{ f\in\mathbb{C}[\mathfrak{g}]:g\cdot f=f\text{ for all }g\in G \right\} .\]
The coinvariant ring of $G$ is the quotient
\[\mathbb{C}[\mathfrak{g}]_G:=\mathbb{C}[\mathfrak{g}]/\mathbb{C}[\mathfrak{g}]_+^G\]
of the full polynomial ring by the ideal $\mathbb{C}[\mathfrak{g}]_+^G$ of $G-$invariant polynomials without constant term. 

Let $\partial_i=\frac{\partial}{\partial{x_i}}$, and for $f\in\mathbb{C}[\mathfrak{g}]$, define $f(\partial)=f(\partial_1,\ldots,\partial_{n^2})$. The $G$-invariant differential operators with no constant term is the set 
\[\mathcal{D}(\mathfrak{g})_+^G=\left\{ f(\partial)|f\in\mathbb{C}[\mathfrak{g}]_+^G \right\}.\]
The module of $G-$harmonic polynomials on $\mathfrak{g}$ is then defined by
\[ \mathcal{H}\left(\mathfrak{g}\right)=\left\{f\in\mathbb{C}[\mathfrak{g}]|\Delta(f)=0\text{ for all }\Delta\in\mathcal{D}(\mathfrak{g})_+^G\right\} .\]
As a (infinite dimensional) representation of $G$, $\mathcal{H}\left(\mathfrak{g}\right)$ is isomorphic to the coinvariant ring $\mathbb{C}[\mathfrak{g}]_G$. 

The $G-$harmonic polynomials also form a graded representation of $GL_n(\mathbb{C})$;
\[ \mathcal{H}\left(\mathfrak{g}\right)=\bigoplus_{d=0}^{\infty}\mathcal{H}\left(\mathfrak{g}\right)^d, \]
where $\mathcal{H}^d\left(\mathfrak{g}\right)=\mathcal{H}\left(\mathfrak{g}\right)\cap\mathbb{C}[\mathfrak{g}]_d$. This fact is non-trivial, since it says that if a function is harmonic, then so are its homogeneous components. In general, every polynomial function can be expressed as a sum of $G-$invariant functions multiplied by $G$-harmonic polynomials. In other words, there is a surjection
\[\mathbb{C}[\mathfrak{g}]^G\otimes\mathcal{H}(\mathfrak{g})\longrightarrow\mathbb{C}[\mathfrak{g}]\longrightarrow0\]
obtained by linearly extending multiplication. Kostant showed in \cite{Kostant:1} that the product of invariants and harmonics is unique, and thus
\[ \mathbb{C}[\mathfrak{g}]\cong \mathbb{C}[\mathfrak{g}]^G\otimes\mathcal{H}(\mathfrak{g}). \]

Each degree $d$ homogeneous component, $\mathcal{H}(\mathfrak{g})^d$, is a finite dimensional rational representation of $G$, and therefore decomposes into a finite direct sum of irreducible representations of $G$. We then pose the question, ``given an irreducible rational representation $V^{\alpha}$ of $GL_n(\mathbb{C})$ with highest weight $\alpha$, in what degrees does $V^{\alpha}$ occur in the direct sum decomposition of $\mathcal{H}(\mathfrak{g})^d$?" This question can be further formalized by considering the \emph{graded multiplicity} of the representation $V^{\alpha}$, which is defined as follows. If $V^{\alpha}$ is a finite dimensional irreducible rational representation of $G$ with highest weight $\alpha$, we denote by $\left[ V^{\alpha},\mathcal{H}(\mathfrak{g})^d\right]$ the multiplicity with which $V^{\alpha}$ occurs in the direct sum decomposition of the homogeneous degree $d$ harmonic polynomials, $\mathcal{H}(\mathfrak{g})^d$. The degrees in which and multiplicities with which the representation $V^{\alpha}$ occurs within the harmonic polynomials is summarized in what is called the \emph{graded multiplicity} polynomial of $V^{\alpha}$ in $\mathcal{H}(\mathfrak{g})$, and is defined by
\begin{equation}\label{eqn:gradedmultiplicity}
    m_{\alpha}(q):=\sum_{d=0}^{\infty}\left[V^{\alpha},\mathcal{H}(\mathfrak{g})^d\right] q^d,
\end{equation} where $q$ is a dummy variable. It follows from \cite{Kostant:1} that $m_{\alpha}(q)\neq 0$ if and only if $\alpha\neq0$ and $\alpha_1+\cdots+\alpha_n=0$. We now state our main result.

\begin{theorem}
\label{theorem1}
Let $V^{\alpha}$ be a $GL_n(\mathbb{C})$ irrep of highest weight $\alpha=(\alpha_1\geq\alpha_2\geq\cdots\geq\alpha_n)$ with $\alpha_1+\alpha_2+\cdots+\alpha_n=0$. Let $k\in\mathbb{N}$ be such that $k\geq |\alpha_n|$, and let $\lambda:=\alpha+(k^n)$. Then $\deg m_{\alpha}(q)=g_{nk,n}(\lambda)$, where $g_{nk,n}$ is the generalized Gini index defined in Equation \ref{eqn:gini2}. 
\end{theorem}

It is worth noting that the Gini index $g_{nk,n}(\lambda)$ is independent of the choice of integer $k$, so long as $k\geq|\alpha_n|$. In simple terms, Theorem \ref{theorem1} states that for a given irreducible rational representation $V^{\alpha}$ of $GL_n(\mathbb{C})$, the graded multiplicity polynomial $m_{\alpha}(q)$ is either zero, or of degree $g_{nk,n}(\lambda)$, where $\lambda=\alpha+(k^n)$, and $k$ is any integer greater than $|\alpha_n|$. In other words, the maximal homogeneous degree in which an irreducible representation of $GL_n(\mathbb{C})$ appears in the $GL_n(\mathbb{C})-$harmonics is completely determined by the Gini index. The proof of this result is provided in Section 3.3.

% Subsection 3.3 - Proof of Main result %%%%%%%%%%%%%%%%%%%%%%%%%%%%%%%%%%%%%%

\subsection{Kostka-Foulkes polynomials, and proof of main result}

If $\lambda=(\lambda_1\geq\cdots\geq\lambda_n)$ is a partition of $nk$ with at most $n$ parts, let $\alpha\in\mathbb{Z}^n$ be given by $\alpha:=\lambda-(k^n)$. Then $\alpha_1+\cdots+\alpha_n=0$ and $k\geq|\alpha_n|$. We will show that, in this setting, the degree of $m_{\alpha}(q)$ is $g_{nk,n}(\lambda)$ by relating the graded multiplicities to \emph{Kostka-Foulkes polynomials}.

Given two partitions $\lambda$ and $\mu$ of a positive integer $n$, the Kostka number $K_{\lambda\mu}$ is the number of semistandard Young tableaux of shape $\lambda$ and weight $\mu$. The Kostka numbers can also be defined as the coefficients obtained by expressing the Schur polynomial $s_{\lambda}$ as a linear combination of monomial symmetric functions $m_{\mu}$;
\[s_{\lambda}=\sum_{\mu}K_{\lambda\mu}m_{\mu}.\]
These numbers are significant in representation theory, as $K_{\lambda\mu}$ counts the dimension of the weight space corresponding to $\mu$ in the irreducible representation $V^{\lambda}$ with highest weight $\lambda$. This result is often expressed in terms of Kostant's partition function $\mathcal{P}:\mathbb{Z}^n\longrightarrow \mathbb{Z}$. Let $\epsilon_i\in\mathbb{Z}^n$ denote the tuple with a 1 in the $i$-th position, and $0$'s elsewhere. Define $\epsilon_{ij}=\epsilon_i-\epsilon_j$. The tuples $\epsilon_{ij}$ with $i\neq j$ form the root system of the Lie algebra $\mathfrak{sl}_n(\mathbb{C})$, and the positive roots are generally defined as the $\epsilon_{ij}$ with $i<j$. Kostant's partition function $\mathcal{P}(\alpha)$ counts the number of ways in which a weight $\alpha$ of $\mathfrak{sl}_n(\mathbb{C})$ can be expressed as a sum of positive roots of $\mathfrak{sl}_n(\mathbb{C})$. This function can also be defined as the coefficients in the series
\[ \sum_{\alpha\in\mathbb{Z}^n}\mathcal{P}(\alpha)x_1^{\alpha_1}\cdots x_n^{\alpha_n}=\prod_{1\leq i<j \leq n}\left(1-\frac{x_i}{x_j}\right)^{-1}_, \]
where $x_1,\ldots,x_n$ are indeterminants. The Kostka numbers can then be written in terms of $\mathcal{P}$, yielding Kostant's weight multiplicity formula
\[ K_{\lambda\mu}=\sum_{w\in S_n}\epsilon(w)\mathcal{P}(w(\lambda+\rho)-(\mu+\rho)), \]
where $\rho=(n-1,n-2,\ldots,1,0)$. The $q-$analogue of Kostant's partition function can be defined as 
\[ \sum_{\alpha\in\mathbb{Z}^n}\mathcal{P}_q(\alpha)x_1^{\alpha_1}\cdots x_n^{\alpha_n}=\prod_{1\leq i<j\leq n}\left( 1-q\frac{x_i}{x_j} \right)^{-1}_. \]
Substituting $\mathcal{P}$ for its $q-$analogue, we obtain the Kostka-Foulkes polynomial
\[K_{\lambda\mu}(q):=\sum_{w\in S_n}\epsilon(w)\mathcal{P}_q(w(\lambda+\rho)-(\mu+\rho)).\]

Hesselink connected $K_{\lambda \mu}(q)$ to the graded multiplicity of an irrep $V^{\alpha}$ of $GL_n(\mathbb{C})$ \cite{Hesselink:1}. If $\alpha\in\mathbb{Z}^n$ is decreasing and sums to zero, and if $\lambda=\alpha+(k^n)$, with $k\geq |\alpha_n|$, then the graded multiplicity polynomial $m_{\alpha}(q)$ defined in Equation \ref{eqn:gradedmultiplicity} is given by
\[ m_{\alpha}(q)=\sum_{w\in S_n}\epsilon(w)\mathcal{P}_q(w(\alpha+\rho)-\rho). \]
Since $\lambda=\alpha+(k^n)$, it follows that
\[m_{\alpha}(q)=K_{\lambda (k^n)}(q).\]
For more on this result see \cite{Kostka-Foulkes:1} and \cite{Gupta:1}.

%The graded multiplicity polynomial $m_{\alpha}(q)$ has a similar formulation due to hesselink.... 
The Kostka-Foulkes polynomials are understood combinatorially in \cite{Lascoux-Schutzenberger:1} and \cite{Macdonald:1}.  From this interpretation, the degree of the polynomial $K_{\lambda \mu}(q)$ can be expressed in terms of the weighted totals of $\lambda$ and $\mu$, defined in Equation \ref{eqn:weightedsum}, and is given by $b(\mu)-b(\lambda)$.

In summary, the only irreducible rational representations $V^{\alpha}$ with highest weight $\alpha\in\mathbb{Z}^n$ for which $m_{\alpha}(q)\neq 0$ are those whose highest weight satisfies $\alpha_1+\cdots+\alpha_n=0$. If $m_{\alpha}(q)$ is nonzero, then we have
\[m_{\alpha}(q)=K_{\lambda(k^n)}(q),\]
where $k$ is any integer satisfying $k\geq|\alpha_n|$ and $\lambda=\alpha+(k^n)$. From \cite{Lascoux-Schutzenberger:1} and \cite{Macdonald:1} It follows that
\[ \text{deg }m_{\alpha}(q)=b((k^n))-b(\lambda)=g_{nk,n}(\lambda). \]
This proves Theorem \ref{theorem1}.

%%%%%%%%%%%%%%%%%%%%%%%%%%%%%%%%%%%%%%%%%%%%%%%%%%%%%%%%%%%%%%%%%%%%%%%%%%%%%%%%%%%%%
% Section 4 - Final Thoughts
%%%%%%%%%%%%%%%%%%%%%%%%%%%%%%%%%%%%%%%%%%%%%%%%%%%%%%%%%%%%%%%%%%%%%%%%%%%%%%%%%%%%%

\section{The Gini index and the Earth Mover's Distance}

\begin{figure}[H]
    \centering
    \includegraphics[width=11cm]{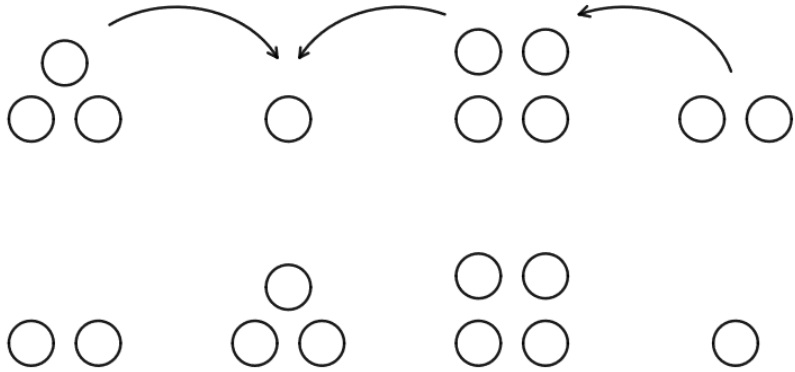}
    \caption{The minimum number of moves required to convert the top distribution, $(3,1,4,2)$, into the bottom distribution, $(2,3,4,1)$, is 3.}
    \label{figure3}
\end{figure}

This formulation of a discrete Gini index was motivated by recent research on the discrete Earth Mover's Distance (EMD) -- also known as the Wasserstein distance. Consider two different distributions of $s$ pebbles into $n$ distinct piles. Simply put, the one-dimensional EMD counts the least number of moves needed to convert the first distribution into the second, where each \emph{move} involves moving one pebble to a neighboring pile. An example of this procedure is provided in Figure \ref{figure3}.
The Gini index and EMD are both measures of dissimilarity and are, in fact, measuring the same quantity.

A composition of a positive integer $s$ into $n$ parts is an $n$-tuple of non-negative integers whose sum is $s$. The EMD of such compositions was discussed in \cite{Willenbring:1} (the one-dimensional case) and \cite{Erickson:2} (a higher dimensional generalization) -- generating functions and connections to representation theory were found in both instances. In \cite{Erickson:1} it was shown that if $\mu$ and $\lambda$ are compositions of $s$ into $n$ parts, then $EMD(\mu,\lambda)$ can be expressed as the symmetric difference of the Young diagrams of the words of $\mu$ and $\lambda$. The word of a composition $\mu$ is a tuple of integers constructed by writing $``0"$ $\mu_1$ times, then $``1"$ $\mu_2$ times, and so on, ending with $``n-1"$ written $\mu_{n}$ times. If $\mu$ is a composition of $s$ into $n$ parts, then the word of $\mu$ will be weakly increasing with length $s$. The Young diagram of a word is a collection of square boxes arranged in left-justified rows, whose ascending row lengths are given by the word. The symmetric difference of the diagrams is the number of boxes in their union minus their intersection.

\begin{example}
\label{example3}
The EMD of the compositions $\mu=(3,1,4,2)$ and $\lambda=(2,3,4,1)$ can by found as in Figure \ref{figure3}, or by using Young diagrams as shown in \cite{Erickson:1}. The word of $\mu$ is $(0,0,0,1,2,2,2,2,3,3)$ and the word of $\lambda$ is $(0,0,1,1,1,2,2,2,2,3)$. Their Young diagrams, and symmetric difference are shown in Figure~\ref{figure4}.
\begin{figure}[H]
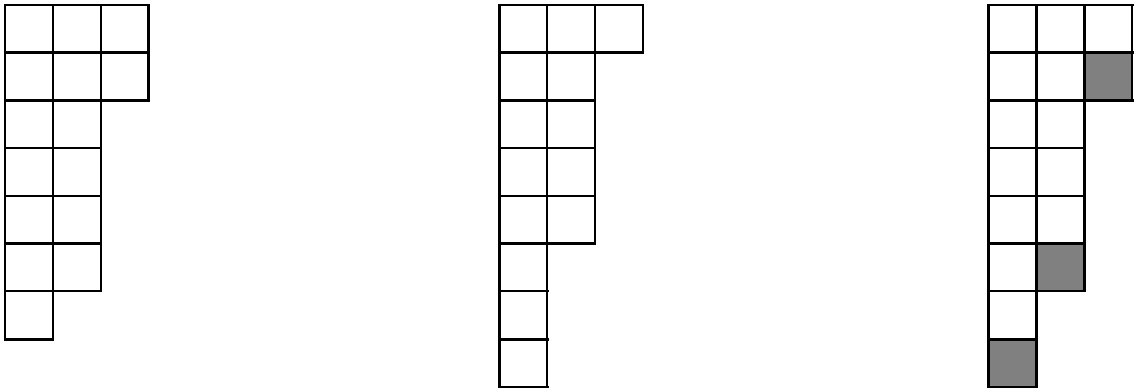

    \centering

    \begin{multicols}{3}
    \ytableausetup{nosmalltableaux}
    \begin{ytableau}
    *(white)  & *(white) & *(white) \\
    *(white)  & *(white) & *(white) \\
    *(white)  & *(white)  \\
    *(white)  & *(white)  \\
    *(white)  & *(white)  \\
    *(white)  & *(white)  \\
    *(white) 
    \end{ytableau}
    \columnbreak
    \ytableausetup{nosmalltableaux}
    \begin{ytableau}
     *(white) & *(white) & *(white) \\
     *(white) & *(white) \\
     *(white) & *(white) \\
     *(white) & *(white) \\
     *(white) & *(white) \\
     *(white) \\
     *(white) \\
     *(white)
    \end{ytableau}
    \columnbreak
    \ytableausetup{nosmalltableaux}
    \begin{ytableau}
     *(white) & *(white) & *(white)  \\
     *(white) & *(white) & *(gray) \\
     *(white) & *(white) \\
     *(white) & *(white) \\
     *(white) & *(white) \\
     *(white) & *(gray) \\
     *(white) \\
     *(gray)
    \end{ytableau}
    \end{multicols}
    
    \caption{The Young diagrams of the words of $\mu$ (left) and $\lambda$ (middle) and their union (right), with symmetric difference shaded.}
    \label{figure4}
\end{figure}
The number of shaded cells in the final diagram is the symmetric difference. Hence, for Example \ref{example3}, $EMD(\mu,\lambda)=3$.
\end{example}
In general, $EMD(\mu,\lambda)$ can be expressed in terms of the weighted totals of $\mu$ and $\lambda$ (defined in Equation \ref{eqn:weightedsum}) as 
\[ EMD(\mu,\lambda)=b(\mu)-b(\lambda) \]
whenever $\mu$ and $\lambda$ are partitions of $s$ into $n$ parts, and $\mu$ majorizes $\lambda$ -- that is, whenever $\sum_{i=1}^k\mu_i \geq \sum_{i=1}^k\lambda_i$ for all $1\leq k \leq n$. The discrete Gini index can also be defined in terms of Young diagrams as in \cite{Kopitzke:1} and \cite{Kopitzke:2}. From this interpretation, in the special case where $s=nk$, it follows that
\[ EMD((k^n),\lambda)=b((k^n))-b(\lambda)=g_{nk,n}(\lambda). \]
That is, the discrete Gini index is a restriction of the discrete one-dimensional Earth Mover's Distance to the set of partitions of $nk$ with at most $n$ parts.

\bibliographystyle{Bibtex.sty}
\bibliography{References}

\bigskip
\hrule
\bigskip

\noindent 2020 {\it Mathematics Subject Classification}:
Primary 05E10, Secondary 91B82, 20G05.

\noindent \emph{Keywords: }
Gini index, Lorenz curve, integer partition, integer composition, general linear group, representation theory, harmonic polynomials, Kostka-Foulkes polynomials, earth movers distance, Wasserstein metric.

\bigskip
\hrule
\bigskip

\end{document}